\documentclass{article}

\usepackage{amsmath}
\usepackage{amssymb}
\usepackage[dvipdfm]{graphicx,color}
\usepackage[all]{xy} 
\usepackage{amsthm, bm}

\newtheorem{theorem}{Theorem}

\newtheorem{remark}{Remark}

\def\phi{\varphi}

\title{Non-integrability criterion for  homogeneous Hamiltonian systems via blowing-up technique of singularities}
\author{Mitsuru Shibayama}

\begin{document}
\maketitle


\begin{abstract}

It is a big problem to distinguish between integrable and non-integrable  Hamiltonian systems.  
We provide a new approach to prove  the non-integrability of homogeneous Hamiltonian systems 
with two degrees of freedom. 
 The homogeneous degree can be chosen from  real values (not necessarily integer).
The proof is based on the blowing-up theory
which McGehee established  in the collinear three-body problem.
We also compare our result with Molares-Ramis theory which is the strongest theory in this field.

\end{abstract}

\section{INTRODUCTION}

Let $H: \mathcal{D}\to \mathbb{R}$ be a smooth function where $\mathcal{D}$ is an open set { in} $\mathbb{R}^{2k}$. 
The Hamiltonian system is defined by the ordinary differential equations 
\begin{equation} \label{eqn:ham}
\frac{dq_{j}}{dt}=\frac{\partial H}{\partial p_{j}}(\mathbf{p}, \mathbf{q}), \quad \frac{dp_{j}}{dt}=-\frac{\partial H}{\partial q_{j}}(\mathbf{p},\mathbf{q}) \qquad (j=1, \dots, k)
\end{equation}
where $(\mathbf{p}, \mathbf{q})=(p_{1}, \dots, p_{k} ,q_{1}, \dots, q_{k}) \in \mathcal{D}$.
 The function $H$ is called the Hamiltonian and 
 $k$ is called the degrees of freedom.

 A function  $F: \mathcal{D} \to \mathbb{R}$ is called the first integral  of \eqref{eqn:ham} if   $F$ is 
 conserved along each solution   of \eqref{eqn:ham}. 
 For two functions $F, G: \mathcal{D} \to \mathbb{R}$, 
the Poisson bracket is the function defined by 
\[ \{ F,  G \}=\sum_{j=1}^{k} \frac{\partial F}{\partial q_{j}}\frac{\partial G}{\partial p_{j}}-\frac{\partial F}{\partial p_{j}}\frac{\partial G}{\partial q_{j}}. \]
{ A function} $F: \mathcal{D} \to \mathbb{R}$ is a first integral
{ of \eqref{eqn:ham}} if and only if 
$\{F, H\}$ is identically zero.
 Hamiltonian system \eqref{eqn:ham} is called integrable if there are $k$ first integrals 
$F_{1}(=H), F_{2}, \dots, F_{k}$ such that 
$dF_{1}, \dots, dF_{k}$ are linearly independent in an open dense set of $\mathcal{D}$
and that $\{ F_{i}, F_{j}\}=0$ for any $i, j = 1, \dots, k$.

The dynamics of the integrable systems are well understood because of 
the Liouville-Arnold theorem(see \cite[Chapter 10]{Arnold}) while
the dynamics of the non-integrable Hamiltonian systems may be chaotic.
Therefore it is important to distinguish between integrable and non-integrable Hamiltonian systems.

This problem have been studied for quite long time.
Bruns \cite{Bruns87} proved that in the 3-body problem there is no algebraic first integral 
which is independent from the known ones.
After that, Poincar\'e \cite{Poincare90} proved that   the
perturbed Hamiltonian systems there is no
 analytic first integral depending analytically on a parameter.
Then by applying it to the restricted 3-body problem, he proved 
the non-existence of an analytic first integral depending analytically on a mass parameter.
 
Another theory in this field was originated by Kovalevskaya \cite{Kovalevskaya89}.
By studying the property of  singularities she  discovered  a new integrable case in the  rigid body model.
As a development of her approach, Ziglin \cite{Ziglin82a, Ziglin82b} established  the theory of singularity for proving the non-integrability.
 By applying the Ziglin analysis, Yoshida \cite{Yoshida87} provided a criterion for the non-integrability of the homogeneous Hamiltonian systems.
Morales-Ruiz \& Ramis \cite{Morales99, MoralesR01}  extended  the Ziglin analysis by applying the Differential Galois theory (Picard-Vessiot theory).
The Morales-Ramis theory is the strongest in this field now.

Our { purpose is 
to prove the non-integrability of  Hamiltonian systems 
from a new approach}.
We consider a Hamiltonian system of 2 degrees of freedom  with   a homogeneous potential of degree $\beta \in \mathbb{R}$.
Its Hamiltonian is { represented by}
\begin{equation} \label{ham}
H(\mathbf{p}, \mathbf{q})=\frac{1}{2} \|\mathbf{p}\|^{2}+U(\mathbf{q}) \qquad ((\mathbf{p}, \mathbf{q}) \in \mathbb{R}^{2} \times { (\mathbb{R}^{2}
\backslash \{\mathbf{0}\})}).  
\end{equation}
Here  $U$ is a real-meromorphic function on $\mathbb{R}^{2} \backslash \{\mathbf{0}\}$ and 
satisfies the homogeneous property:
\[ U(\lambda \mathbf{q})=\lambda^{\beta} U(\mathbf{q})  \qquad (\mathbf{q} \in \mathbb{R}^{2} \backslash \{\mathbf{0}\},\lambda >0).\]

Let $V(\theta)=U(\cos \theta, \sin \theta)$.

\begin{theorem}\label{th1}
Assume the following 6 properties:
\begin{enumerate}
\item the homogeneous degree $\beta$ is a real number excluding $-2$ and $0$: 
\[ \beta \in \mathbb{R}\backslash \{-2, 0\};\]
  \item there are three critical points  $\theta_{l}$ of $V$: 
  \[ \frac{\partial V}{\partial \theta}(\theta_{l})=0, \quad \theta_{-1}<\theta_{0}<\theta_{1}, \quad ( l=-1,0,1);\]
  \item the function $V$ is negative between $\theta_{-1}$ and $\theta_{1}$: 
  \[  V(\theta) < 0    \quad ( \theta \in [\theta_{-1},\theta_{1}]);\] 
  \item the derivative of $V$ does not vanish between these critical points:
  \[ \frac{\partial V}{\partial \theta}(\theta ) \neq 0 \quad (\theta \in (\theta_{-1},\theta_{0}) \cup (\theta_{0}, \theta_{1})); \]
  \item the second derivative of $V$ is negative at critical points $\theta_{\pm 1}$: 
  \[\frac{\partial^2 V}{\partial \theta^2}(\theta_{\pm1})<0;\]
  \item at critical point $\theta_{0}$, the following inequality satisfies: 
\[-\frac{1}{8}(\beta+2)^{2}V(\theta_{0})<\frac{\partial^2 V}{\partial \theta^2}(\theta_{0}).\]
\end{enumerate}
Then the Hamiltonian system of \eqref{ham} has no real-meromorphic  first integral independent from $H$. 

\end{theorem}

\begin{figure}[htbp]
\begin{center} 
\includegraphics[width=8cm]{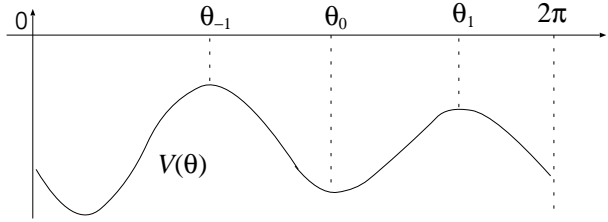}
\caption{Function $V(\theta)$}
\end{center}
\end{figure}

Above we used the word ``real-meromorphic''.
We call a real function $f(\mathbf{p}, \mathbf{q})$
real-meromorphic if 
and only if 
$f(\mathbf{p}, \mathbf{q})$ 
is analytic in all but possibly a discrete subset of $\mathbb{R}^{2} \times 
(\mathbb{R}^{2} \backslash \{\mathbf{0}\})$
and these exceptional points must be poles.

\begin{remark}
The case of $\theta_{1}=\theta_{-1}+2\pi$ 
is allowed in assumption 2.
These two critical points are essentially identical.
In this case, just two critical points of $V$ are necessary.
\end{remark}

\begin{remark}\label{rem:-2}
In the case of  $\beta=-2$, the Hamiltonian system is  integrable.
Because a function 
\[ G(\mathbf{p}, \mathbf{q})=(\mathbf{q} \cdot \mathbf{p})^{2} -2 \|\mathbf{q}\|^{2}H(\mathbf{p}, \mathbf{q})\]
is a first integral.
Hence this case does not need to be studied.
\end{remark}

\begin{remark}\label{rem:pos}

In the case of $V(\theta) > 0$ on $[\theta_{-1}, \theta_{1}]$,
if $V$ is analytic in the complex domain $\mathbb{C}^{2} \backslash \{(0, 0)\}$, 
$V$ can be replaced  by changing coordinates with  $( \mathbf{P}, \mathbf{Q})=(\sqrt{-1} \mathbf{p}, \sqrt{-1}\mathbf{q})$,
and then the new equations satisfy the assumption 2 of this theorem.
\end{remark}

If $V$ is a constant, 
the system is integrable.
Hence we need to consider the non-constant functions. 
Generically there are several critical points of $V$
and the graph is convex at some of them.
The assumption 1-5 of this theorem is not  strong, and 
only  assumption 6   is a little strong.

This paper is organized as follows. 
 In Section 2 we introduce the McGehee's blowing-up technique for  the homogeneous Hamiltonian systems.
We prove  our theorem in Section 3 by using the McGehee's technique.
We present two applications of the theorem in Section 4. 
In the final section we compare our theorem with the Morales-Ramis theorem.

\section{MCGEHEE'S BLOWING UP TECHNIQUE}

McGehee \cite{McGehee74} established a  blowing-up technique for the triple collision singularity in the 
collinear three-body problem. 
We can easily extend the technique for the general homogeneous Hamiltonian systems \eqref{ham}. 

We fist consider the case of $\beta <0$.
The McGehee  coordinates $(r, \theta, v, w)$  are defined by 
\[ \mathbf{q}=r(\cos \theta, \sin \theta), \mathbf{p}=r^{\beta/2}(v(\cos \theta, \sin \theta)+w(-\sin \theta, \cos \theta))\]
and 
the time variable $t$ is changed into   $\tau$ according to $dt = r^{1-\beta/2} d \tau$. 
The map  $(r, \theta, v, w) \mapsto (\mathbf{p}, \mathbf{q})$ are analytic 
{ in $\{(r, \theta, v, w) \mid r >0, \theta \in \mathbb{R}/2 \pi \mathbb{Z}, v, w \in \mathbb{R}\}$}.
Then the equations become
\begin{align}
\frac{d r}{d \tau}&=r v  \label{dr} \\
\frac{d\theta}{d \tau}& = w \label{dth} \\
 \frac{d v}{d \tau}  &=-\frac{\beta}{2}  v^{2}+  w^{2} -  \beta  V(\theta) \label{dv}\\
  \frac{d w}{d \tau}  &=-\left(\frac{\beta}{2}+1\right) vw -  \frac{\partial V}{\partial \theta}(\theta). \label{dw}
\end{align}
In these coordinates the total energy is
\begin{equation}\label{energy}
h=r^{\beta} \left(\frac{v^{2}+w^{2}}{2} +V(\theta)\right).
\end{equation}
Fix the energy constant at any non-zero value($h \neq 0$).

The point $\mathbf{q}=\mathbf{0}$ is singularity of the differential equations,
 but $r =0$ is not singular in these differential equations \eqref{dr}-\eqref{dw}.
It is sufficient to consider the three equations 
 \eqref{dth}, \eqref{dv} and \eqref{dw}, 
 since these equations are independent from $r$ and 
 since $r$ can be obtained from 
 \eqref{energy}.

The  set 
\[ \mathcal{M}=\left\{ ( \theta, v, w) \in \mathbb{R}/ 2 \pi \mathbb{Z} \times \mathbb{R} \times \mathbb{R} \mid  \frac{v^{2}+w^{2}}{2} +V(\theta)=0\right\} \]
 is invariant.
 In the case of the $n$-body problem, 
 $\mathcal{M}$ is 
  called the collision manifold.
Orbits converge to  $\mathcal{M}$ as $r \to 0$.

In the case that $\beta >0$, 
we can discuss similar argument by 
letting $R=r^{-1}$.
The equation \eqref{dr} becomes 
\begin{equation}\label{dR}
 \dot{R}=-Rv
 \end{equation}
and the total energy is 
\begin{equation}\label{energy2}
 H=R^{-\beta} \left(\frac{v^{2}+w^{2}}{2} +V(\theta)\right).
 \end{equation}
The equations can be extended to $R=0$.
Orbits converge to the invariant set $\mathcal{M}$ as $R \to 0$.
It is sufficient to consider the three equations 
 \eqref{dth}, \eqref{dv} and \eqref{dw}.

The flow on $\mathcal{M}$ is gradient-like if $\beta \neq -2$.
This means that 
the $v$-component is monotone along each solution  excluding  equilibrium points 
since all orbits on $\mathcal{M}$  satisfy
\[ \frac{ d v}{d \tau}=\left(\frac{\beta}{2}+1\right)w^{2} \begin{cases} 
\ge 0 & (\beta > -2) \\
\le 0 & (\beta < -2)
\end{cases}  \] 

If  $\theta_{c}$ is a critical point of $V$, i.e. $\frac{\partial V}{\partial \theta}(\theta_{c})=0$,
 $(\theta, v, w)=(\theta_{c}, \pm \sqrt{-2V(\theta_{c})}, 0)$ are equilibrium points of \eqref{dth}, \eqref{dv}, \eqref{dw}.
The linearized equations of \eqref{dth}, \eqref{dv}, \eqref{dw}  at $(\theta, v, w)=(\theta_{c},  \pm \sqrt{{-}2V(\theta_{c})}, 0)$ are
  \begin{equation*}
  \frac{d}{d \tau}\left(\begin{array}{c}\delta \theta \\\delta v \\\delta w\end{array}\right)=
  \left(\begin{array}{ccc}0 & 0 & 1 \\0 &  \mp \beta\sqrt{{-}2V(\theta_c)} & 0 \\ {-}\frac{\partial^{2}V}{\partial \theta^{2}}(\theta_c) & 0 & \mp \left(\frac{\beta}{2}+1\right) \sqrt{{-}2V(\theta_{c})} \end{array}\right)
  \left(\begin{array}{c}\delta \theta \\\delta v \\\delta w\end{array}\right).
  \end{equation*}
The eigenvalues of the coefficient matrix are 
$\lambda_{1}=\mp \beta\sqrt{2V(\theta_c)}$, $\lambda_{2}$ and $\lambda_{3}$ 
where $\lambda_{2}$ and $\lambda_{3}$ are the roots of equation
\[ \lambda^{2} \pm \left(\frac{\beta}{2}+1\right) \sqrt{{-}2V(\theta_{c})} \lambda{+}\frac{\partial^{2}V}{\partial \theta^{2}}(\theta_c)=0.\]
The eigenspace corresponding to $\lambda_{1}$ is perpendicular to
$\mathcal{M}$ at the equilibrium point
and the eigenspace corresponding to $\lambda_{2}$ and $\lambda_{3}$ 
is tangent to $\mathcal{M}$.

\section{PROOF OF THEOREM 1}

Assume that $\Phi(\mathbf{p}, \mathbf{q})$ is a real-meromorphic first integral of \eqref{ham}.
From the homogeneous property
if $(\mathbf{p}(t), \mathbf{q}(t))$ is a solution, so is $(c^{\beta}\mathbf{p}(c^{\beta-2}t), c^{2} \mathbf{q}(c^{\beta-2} t))$ 
for any constant $c>0$.
Then $\Phi(c^{\beta} \mathbf{p}, c^{2} \mathbf{q})$ is also an first integral.

The point $(\mathbf{p}, \mathbf{q})=(\mathbf{0}, \mathbf{0})$ may be an essential singularity of $\Phi$.
Consider the Laurent series at this point:
\[ \Phi(\mathbf{p}, \mathbf{q})=\sum_{k_{1}, k_{2}, k_{3}, k_{4} \in \mathbb{Z}}
a_{k_{1}k_{2}k_{3}k_{4}} p_{1}^{k_{1}}p_{2}^{k_{2}}q_{1}^{k_{3}}q_{2}^{k_{4}}. \]
Then we get 
\[ \Phi(c^{\beta} \mathbf{p}, c^{2} \mathbf{q})=\sum_{k_{1}, k_{2}, k_{3}, k_{4} \in \mathbb{Z}}
a_{k_{1}k_{2}k_{3}k_{4}} c^{\beta(k_{1}+k_{2})+2(k_{3}+k_{4})} p_{1}^{k_{1}}p_{2}^{k_{2}}q_{1}^{k_{3}}q_{2}^{k_{4}}. \]
We gather the terms according to the power of $c$
\begin{equation} \label{Lau}
 \Phi(c^{\beta} \mathbf{p}, c^{2} \mathbf{q})=\sum_{\omega \in \Omega} c^{\omega} f_{\omega} (\mathbf{p}, \mathbf{q})
 \end{equation}
 where 
 \[ \Omega=\{ \beta(k_{1}+k_{2})+2(k_{3}+k_{4}) \mid k_{j} \in \mathbb{Z}, a_{k_{1}k_{2}k_{3}k_{4}}\neq 0 \}\]
 and
 \[ f_{\omega}(\mathbf{p}, \mathbf{q})=
 \sum_{\beta(k_{1}+k_{2})+2(k_{3}+k_{4}) =\omega} a_{k_{1}k_{2}k_{3}k_{4}} p_{1}^{k_{1}}p_{2}^{k_{2}}q_{1}^{k_{3}}q_{2}^{k_{4}}.\]
 By substituting $bc$ for $c$ of \eqref{Lau}, we get
\begin{equation} \label{eqn:bc1}
\Phi(b^{\beta}c^{\beta}  \mathbf{p}, b^{2}c^{2} \mathbf{q})=\sum_{\omega \in \Omega} b^{\omega} c^{\omega} f_{\omega} (\mathbf{p}, \mathbf{q}), 
\end{equation}
and by substituting $c$,
$\mathbf{p}$ and $\mathbf{q}$ for $b$, $c^{\beta} \mathbf{p}, c^{2} \mathbf{q}$ of \eqref{Lau}, we get
\begin{equation}\label{eqn:bc2}
 \Phi(b^{\beta}c^{\beta}  \mathbf{p}, b^{2}c^{2} \mathbf{q})=\sum_{\omega \in \Omega} b^{\omega} f_{\omega} (c^{\beta}\mathbf{p}, c^{2}\mathbf{q}).
 \end{equation}
These equations \eqref{eqn:bc1} and \eqref{eqn:bc2} deduce  
\[ \sum_{\omega \in \Omega} b^{\omega} f_{\omega} (c^{\beta}\mathbf{p}, c^{2}\mathbf{q})=\sum_{\omega \in \Omega} b^{\omega} c^{\omega} f_{\omega} (\mathbf{p}, \mathbf{q}).\]
Therefore we get 
\[ f_{\omega}(c^{\beta} \mathbf{p}, c^{2} \mathbf{q})=c^{\omega}f_{\omega} (\mathbf{p}, \mathbf{q}).\]
Moreover since
\[ \frac{d }{ dt }\Phi(c^{\beta} \mathbf{p}, c^{2} \mathbf{q})=\sum_{\omega \in \Omega}
c^{\omega}\frac{d}{dt} f_{\omega} (\mathbf{p}, \mathbf{q}) =0\] 
for any $c$, each
$f_{\omega} (\mathbf{p}, \mathbf{q})$ is a first integral.

Therefore we can assume that the first integral $\Phi$ satisfies
\begin{equation}\label{finthom}
\Phi(c^{\beta} \mathbf{p}, c^{2} \mathbf{q})=c^{\rho} \Phi(\mathbf{p}, \mathbf{q})
\end{equation}
for some constant $\rho$.

From here we focus the case of $-2 < \beta <0$.
Let 
\[ \Psi (r, \theta, v, w)=\Phi(r^{-\beta/2}(v\cos \theta -w \sin \theta), r^{-\beta/2}(v \sin \theta +w \cos \theta), r\cos \theta, r \sin \theta). \]
From the property \eqref{finthom},  
$\Psi$ can be written by 
\[ \Psi(r, \theta, v, w)=r^{\rho/2} \Psi(1, \theta, v, w).\]
The function $\Psi(1, \theta, v, w)$ is real-meromorphic
of $(\theta, v , w)$.
Note that we do not need analyticity at $r=0$ because of $r=1$.

We denote the equilibrium points by
\[ D_{l}^{\pm} =(\theta_{l}, \pm \sqrt{-2 V(\theta_{l})}, 0) \qquad (l=-1, 0,1).\]
We also use local coordinates $(\theta, w, z)$ near $D_{l}^{-}$ where
\[  z=\frac{v^{2}+w^{2}}{2}+V(\theta). \]
The transformation 
$\{(\theta, v, w) \mid \theta \in \mathbb{R}/2\pi \mathbb{Z}, v <0, w \in \mathbb{R} \} 
\to \{(\theta, z, w) \mid \theta \in  \mathbb{R}/2\pi \mathbb{Z}, z \ge \frac{w^{2}}{2}+V(\theta)  \} $
is real-analytic.
The surface $\mathcal{M}$ corresponds to the plane $z=0$.
In these coordinates, 
the energy is  represented by
\[ h=r^{\beta}z.\]
Define a function $g$ on a neighborhood by 
\[ g(\theta, z, w)=\Psi(1, \theta, -\sqrt{2z-w^{2}-2V(\theta)}, w)\]
which is real-{meromorphic} where the coordinates work.
Because $\Psi$ is real-meromorphic, we can consider the Laurent series of $g$ at $z=0$ with respect to 
$z$: 
\[ g=\sum_{k=\nu}^{\infty} \gamma_{k}(\theta, w) z^{k}\]
where $\nu$ is an integer and  $\gamma_{\nu}(\theta, w)$ is not identically zero.
Hence the first integral is represented by
\begin{align*}
 \Psi(\big(\frac{h}{z}\big)^{\frac{\rho}{2\beta}}, \theta, -\sqrt{2z-w^{2}-2V(\theta)}, w)&=\big(\frac{h}{z}\big)^{\frac{\rho}{2\beta}}\sum_{k=\nu}^{\infty} \gamma_{k}(\theta, w) z^{k} =:\Xi(\theta, w, z). 
 \end{align*}
 If $\Phi$ depends only on $H$, $\Xi$ is a constant function.
  From here the proof varies according to $\nu-\frac{\rho}{2 \beta}$.

\paragraph{The case of $\nu-\frac{\rho}{2 \beta} <0$.}

\begin{figure}[htbp]
\begin{center}
\includegraphics[width=8cm]{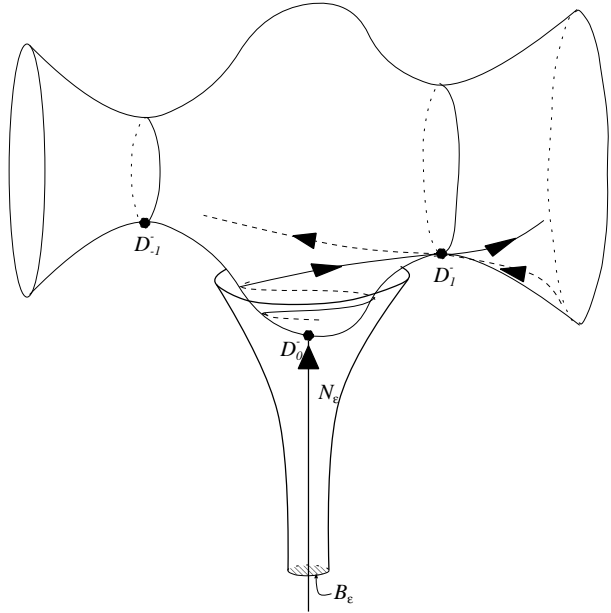}
\caption{A  solution converging to $D_{0}^{-}$}
\end{center}
\end{figure}

Take any $P \in W^{s}(D_{0}^{-}) \backslash \mathcal{M}$ near $D_{0}^{-}$.
Let $a=\Xi(P)$.
We take a small neighborhood of $P$
\[ B_{\varepsilon}=\{Q \in \mathbb{R}^{3}  \mid |P-Q| <\varepsilon \} \qquad (0 < \varepsilon \ll 1), \]
such that  for any $Q \in B_{\varepsilon}$, 
\begin{equation}\label{ene}
 a-1 \le \Xi(Q) \le a+1
 \end{equation}
is satisfied.
Let $\phi_{\tau}(\theta, z, w)$ be the flow of the differential equations.
Since the first integral is conserved along each orbit, 
\eqref{ene} holds in
\[ N_{\varepsilon}=\{ \phi_{\tau} (Q) \mid \tau \ge 0, Q \in B_{\varepsilon}\}.\]

From the continuity,  
\eqref{ene} also holds its closure $\overline{N}_{\varepsilon}$.
This set $\overline{N}_{\varepsilon}$ includes the unstable manifold $W^{u}(D_{0}^{-})$ of $D_{0}^{-}$, 
and $W^{u}(D_{0}^{-})$ is an open set of $\mathcal{M}$. 
The $z$-component converges to zero as $Q$ goes close to $\mathcal{M}$.
Hence  $\gamma_{\nu}$ must be zero on $W^{u}(D_{0}^{-})$.
From the analyticity,  $\gamma_{\nu}$ is identically zero.
This contradicts the assumption for $\gamma_{\nu}$.

\paragraph{The case of $\nu-\frac{\rho}{2 \beta} >0$.}

\begin{figure}[htbp]
\begin{center}
\includegraphics[width=8cm]{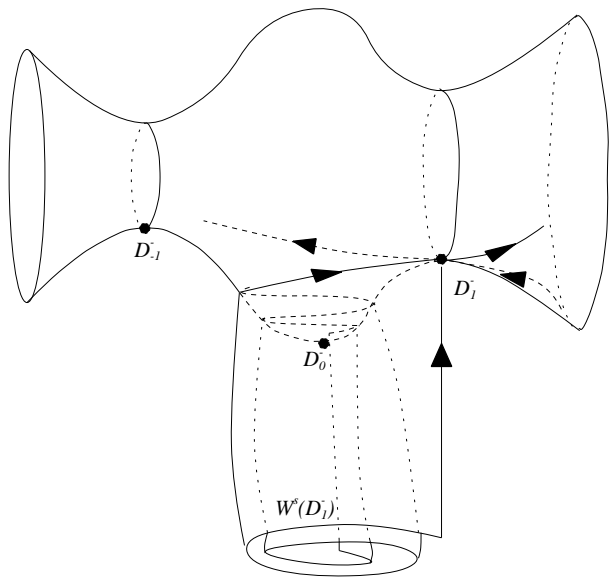}
\caption{The stable manifold of $D_{1}^{-}$}
\end{center}
\end{figure}

Consider the case of  $V(\theta_{1}) \le V(\theta_{-1})$.
The other case is essentially same.
Take any  $Q \in W^{s}(D_{1}^{-}) \backslash \mathcal{M}$.
The first integral has a value $c$ along the orbit passing $Q$:
 \[ \Xi(\varphi_{\tau}(Q))=c \qquad (\tau \in \mathbb{R}).\]
The $z$-component of $\varphi_{\tau}(Q)$ converges to $0$  as $\tau$ diverges to infinity, 
then $c$ must be $0$.
Therefore $\Xi(Q)$ is zero for all $Q \in W^{s}(D_{1}^{-}) \backslash \mathcal{M}$.
The closure of $W^{s}(D_{1}^{-}) \backslash \mathcal{M}$ includes $W^{s}(D_{1}^{-})$.
Because of the continuity, $\Xi(Q)$ is zero on $W^{s}(D_{1}^{-})$.
We can write the function $\Xi$ as 
\[ \Xi(\theta, w, z)=\big(\frac{h}{z}\big)^{\frac{\rho}{2\beta}}z^{\nu}\sum_{k=0}^{\infty} \gamma_{k+\nu}(\theta, w) z^{k}. \]
Therefore 
\[ \sum_{k=0}^{\infty} \gamma_{k+\nu}(\theta, w) z^{k} =\gamma_{\nu}(\theta, w)
+\gamma_{\nu+1}(\theta, w)z+\dots=0\]
satisfies 
on $W^{s}(D_{1}^{-}) \backslash \mathcal{M}$. 
  From the continuity, $\gamma_{\nu}=0$ on $W^{s}(D_{1}^{-}) \cap \mathcal{M}$.

Since $\frac{\partial^2 V}{\partial \theta^2}(\theta_{1}) <0$, the equilibrium point $D_{1}^{-}$ is hyperbolic and $\lambda_{2} \lambda_{3} <0$.
Hence
there are  stable and unstable manifolds with dimension 1  on $\mathcal{M}$.
The dynamics near the equilibrium point $D_{0}^{-}$ is stable focus
and the flow on $\mathcal{M}$ is gradient-like with respect to the $v$-component.
Hence $W^{u}(D_{1}^{-})$ twins around $D_{0}^{-}$
and $\Xi$ is equal to zero on the spiral curve.
$\gamma_{\nu}$ is also zero there.
Therefore from analyticity $\gamma_{\nu}(\theta, w) \equiv 0$. 
This is a contradiction.

\paragraph{The case of $\nu-\frac{\rho}{2 \beta} =0$.}
In this case 
$\gamma_{\nu}$ is a first integral for the flow on $\mathcal{M}$.
From the similar argument as the previous case, 
$\gamma_{\nu}$ is a constant { $c$}.
$\Xi  - c$ is also a first integral.
If $\Xi-c$ is not identically zero, 
$\Xi -c$  has zero point of finite degree  at $z=0$.
This is reduced  to  the  case of $\nu-\frac{\rho}{2 \beta} >0$.
This completes the proof for $-2 < \beta <0$. \\

The proof for the other $\beta$ is essentially same.
We survey the cases.

Consider the case of  $\beta < -2$.
\begin{description}
\item[The case of $\nu-\frac{\rho}{2 \beta} <0$] 
$\gamma_{\nu}$ must be zero $W^{u}(D_{l}^{-})$. 
One branch of $W^{u}(D_{\pm 1}^{-})$ twins around $D_{0}^{-}$. 
\item[The case of  $\nu-\frac{\rho}{2 \beta} >0$] 
$\gamma_{\nu}$ must be zero $W^{s}(D_{l}^{-})$.
Since $W^{s}(D_{0}^{-})$ is an open set of $\mathcal{M}$, 
$\gamma_{\nu}$ must be a zero function.
\item[The case of $\nu-\frac{\rho}{2 \beta} =0$] 
$\gamma_{\nu}$ must be constant  $W^{s/u}(D_{l}^{-})$.
If $\Xi$ is not constant function, this case can be reduced to the case of $\nu-\frac{\rho}{2 \beta}>0$.
\end{description}

Finally consider the case of $\beta >0$.
\begin{description}
\item[The case of $\nu+\frac{\rho}{2 \beta} <0$] 
$\gamma_{\nu}$ must be zero $W^{s}(D_{l}^{-})$.
One branch of $W^{s}(D_{\pm 1}^{-})$ twins around $D_{0}^{-}$. 
\item[The case of $\nu+\frac{\rho}{2 \beta} >0$]
$\gamma_{\nu}$ must be zero $W^{u}(D_{l}^{-})$. 
Sine $W^{u}(D_{0}^{-})$ is an open set of $\mathcal{M}$.
$\gamma_{\nu}$ must be a zero function.
\item[The case of  $\nu+\frac{\rho}{2 \beta} =0$]
$\gamma_{\nu}$ must be constant  $W^{s/u}(D_{l}^{-})$.
If $\Xi$ is not constant function, this case can be reduced to the case of $\nu+\frac{\rho}{2 \beta}>0$.
\end{description}

\section{APPLICATION}
\paragraph{The Isosceles Three-Body Problem}

\begin{figure}[htbp]
\begin{center}
\includegraphics[scale=0.7]{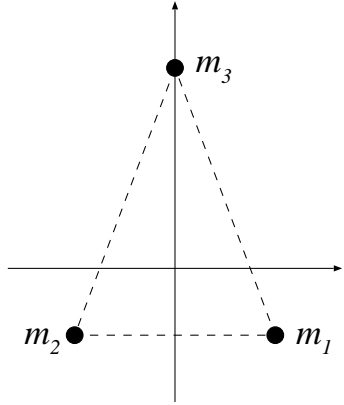}
\caption{The planar isosceles three-body problem}\label{fig:pliso}
\end{center}
\end{figure}

In the planar isosceles three-body problem,
 we can take the centre of gravity  as the origin
 and the symmetric axis as the $y$-axis,
 and assume that the equal masses are located at
\begin{equation*}
(x,y)\quad\mbox{and}\quad(-x,y)
\label{eqn:equal}
\end{equation*}
and the other mass $m_3$ is located at
\begin{equation*}
(0,-2\alpha^{-1}y) \hspace{7cm}
\label{eqn:third}
\end{equation*}
in the inertial coordinate system, where $\alpha=m_3/m$(Figure \ref{fig:pliso}).

 By rescaling it, 
the Hamiltonian is represented by
\[ H(\mathbf{p}, \mathbf{q})=\frac{1}{2} (p_{1}^{2}+p_{2}^{2})-\frac{1}{q_{1}}
 -\frac{4\alpha^{3/2}}{\sqrt{\alpha q_{1}^{2} +(\alpha+2)q_{2}^{2}}}.\]

By applying Theorem \ref{th1}, we obtain:

\begin{theorem}
If  $\alpha < \frac{55}{4}$, 
 the isosceles three-body problem has no real-meromorphic first integral 
independent from $H$.
\end{theorem}

In fact, it is known that 
 the dynamics is complex
in the case of  $\alpha < \frac{55}{4}$.
For example 
there are infinitely many heteroclinic orbits\cite{Moeckel, SY}.

\paragraph{Yoshida's Example}
Consider the Hamiltonian
\begin{equation}\label{Yoshidaex}
 H(\mathbf{p}, \mathbf{q})=\frac{1}{2}(p_{1}^{2}+p_{2}^{2}){+} \frac{1}{4}(q_{1}^{4}+q_{2}^{4})
{ +}\frac{\varepsilon}{2} q_{1}^{2}q_{2}^{2},  
 \end{equation}
which was written  on Yoshida's paper \cite{Yoshida86}. 
{
As we stated at Remark \ref{rem:pos}, 
we can consider the Hamiltonian
\begin{equation}\label{Yoshidaex2}
 G(\mathbf{p}, \mathbf{q})=\frac{1}{2}(p_{1}^{2}+p_{2}^{2})- \frac{1}{4}(q_{1}^{4}+q_{2}^{4})
-\frac{\varepsilon}{2} q_{1}^{2}q_{2}^{2}
 \end{equation}
 instead of $H$.
By applying Theorem \ref{th1}, we obtain:
\begin{theorem}\label{th:3}
If  $\varepsilon < -\frac{1}{8}$ or $\varepsilon >  \frac{25}{7}$, 
the Hamiltonian system \eqref{Yoshidaex2} has  no real-meromorphic first integral independent from $G$. 
\end{theorem}

From Theorem \ref{th:3} and Remark \ref{rem:pos}, we obtain: 
\begin{theorem}
If  $\varepsilon < -\frac{1}{8}$ or $\varepsilon >  \frac{25}{7}$, 
the Hamiltonian system \eqref{Yoshidaex} has  no meromorphic first integral independent from $H$. 
\end{theorem}

\section{COMPARISON WITH THE MORALES-RAMIS THEORY}

We call a  configuration $\mathbf{c} \in \mathbb{R}^{2}$ the Darboux point of $U$ if  
$\nabla U(\mathbf{c}) =\mathbf{c}$.
Consider the Hessian matrix of  $U$ at $\mathbf{c}$ and 
call its eigenvalues Yoshida coefficients at $\mathbf{c}$.
Since $U$ is homogeneous with degree $\beta$, 
we can easily show that one of Yoshida coefficients is $\beta-1$.
As computed by Sansaturio et al \cite{Sansaturio97},
the other  (non-trivial) Yoshida coefficient is  represented by
\[  \lambda=\beta^{-1} V(\theta_{c})^{-1} \frac{\partial^{2}V}{\partial \theta^{2}}(\theta_{c}) +1\]
in the polar coordinates where $\frac{\partial V}{\partial \theta}(\theta_{c})=0$.

In our theorem   the assumption 6 can be written as 
\[ -\frac{1}{8}(\beta+2)^{2} > (\lambda-1) \beta,\]
by using $\lambda$.
Then, in other words, if an integrable Hamiltonian system satisfies the assumption 1-5, 
the Yoshida coefficients at each Darboux point satisfy 
\begin{equation}\label{ourcri}
 -\frac{1}{8}(\beta+2)^{2} \le (\lambda-1) \beta.
 \end{equation}
The Morales-Ramis theorem gave a list of the Yoshida coefficient 
which integral systems can have.
We have compared the inequality \eqref{ourcri} and the Morales-Ramis' list.
The integrable list  given by Morales-Ramis is included in our region \eqref{ourcri} for $\beta \in \mathbb{Z} \backslash \{\pm 2, 0\}$.
For example, in the case of  
$\beta=-1$, 
from the Moreles-Raims theorem, the Yoshida coefficient of an integrable system 
must be in
\[ \{  -\frac{1}{2} p (p-3) \mid p \in \mathbb{Z}\}=\{1, 0, -2, -5, -9, \dots\}. \]
According to our theorem, the Yoshida coefficient of an integrable system 
must be no more than $9/8$ if the other assumptions 1-5 are satisfied.

In the example of the isosceles three-body problem, 
the Morales-Ramis theory guarantees the non-existence of 
meromorphic first integral for any $\alpha$.
In the Yoshida's example, 
Morales-Ramis theory guarantees the non-existence of 
meromorphic first integral excluding $\varepsilon=0,1,3$.
The same result have been obtained through the Ziglin analysis \cite{Yoshida86}.
It is known that these exceptional three cases are actually integrable.

We compare our theorem with the Morales-Ramis theory
in several viewpoints.

\paragraph{Homogeneous degree}
Our theorem can be applied to the case of any {\it real number} $\beta$ excluding $-2, 0$
 while the result from an application \cite{Morales99} of Morales-Ramis theory
can be apply to the case of any {\it integer} excluding $\beta=-2, 0, 2$.
The case of $\beta=-2$ does not need to be studied since the systems are integrable
as we stated at Remark \ref{rem:-2}.
Our theorem alone can be applied to the case of $\beta=2$
\footnote{Recently Andrzej J. Maciejewski provided 
a non-integrability criterion in the case of $\beta=2$ by using the Morales-Ramis theory and 
the higher-order variational equation in his talk at some conference.}.
{ Neither  show anything in the case of $\beta=0$.}

\paragraph{Degrees of freedom}
Our theorem can be applied to {\it two} degrees of freedom while 
Morales-Ramis theory can be applied to {\it any} degrees of freedom.

\paragraph{Yoshida coefficients} 

In the case of integer $\beta$ except $0, \pm 2$, 
the assumption which is imposed in the Morales-Ramis theory is wider than ours 
for proving the non-integrability.

\paragraph{Class of functions}

Our function class of first integrals is bigger. 
We prove the non-existence of first integral which is {\it meromorphic 
as a real function in $\mathbb{R}^{2} \times (\mathbb{R}^{2}\backslash \{(0,0)\})$}, 
while M-R theory 
prove the non-existence of first integrals which is 
{\it meromorphic 
as a complex function}.
Moreover  only our class of functions allows essential singularities at the exceptional points: 
$\mathbf{q}=\mathbf{0}, \mathbf{q}=\infty, \mathbf{p}=\infty$.

\paragraph{Proof methods}  
Proofs are quite different.
Our proof is  simpler and based on dynamics (the behavior of stable and unstable manifolds). 
the proof of Morales-Ramis theory  is  far from the theory of the dynamics
since that is based on the complex analysis and the differential Galois theory.
\\~\\
{\bf Acknowledgement}
The author is  supported by the Sumitomo Foundation, Grant for Basic Science Research Projects 
No. 111153.

\end{document}